\documentclass{article}

\usepackage{arxiv}

\usepackage[utf8]{inputenc} 
\usepackage[T2A]{fontenc}
\usepackage{url}  
\usepackage{booktabs} 
\usepackage{amsfonts} 
\usepackage{nicefrac} 
\usepackage{microtype} 
\usepackage{lipsum}

\usepackage{graphicx}
\usepackage{amsmath}
\usepackage{xcolor}

\usepackage[english]{babel}
\usepackage[ruled,vlined]{algorithm2e}

\usepackage{cite}

\title{Introducing phase jump tracking - a fast method for eigenvalue evaluation of the direct Zakharov-Shabat problem}

\author{
 Igor Chekhovskoy$^{1,2,*}$, Sergey Medvedev$^{2,1}$, Irina Vaseva$^{2,1}$, Egor Sedov$^{1}$, Mikhail Fedoruk$^{1,2}$\\
$^{1}$ Novosibirsk State University, Novosibirsk 630090, Russia,\\
$^{2}$ Institute of Computational Technologies, SB RAS, Novosibirsk
630090, Russia,\\
* Corresponding author: i.s.chekhovskoy@nsu.ru
}

\begin{document}
\maketitle

\begin{abstract}
We propose a new method for finding discrete eigenvalues for the direct Zakharov-Shabat problem, based on moving in the complex plane along the argument jumps of the function~$a(\zeta)$, the localization of which does not require great accuracy.
It allows to find all discrete eigenvalues taking into account their multiplicity faster than matrix methods and contour integrals. The method shows significant advantage over other methods when calculating a large discrete spectrum, both in speed and accuracy.
\end{abstract}

\keywords{Zakharov-Shabat problem \and Direct scattering transform \and Nonlinear Fourier transform \and Nonlinear Schr\"odinger equation \and Fast numerical methods}

\section*{Introduction}

The method of the inverse scattering problem (ISM) for the nonlinear Schr\"odinger equation (NLSE) allows to integrate this equation analytically~\cite{ZakharovShabat1972, Ablowitz1981} and consists of two steps: the first step is solving the direct Zakharov-Shabat problem (ZSP), i.e. we need to find scattering data, and the second step is solving the inverse scattering problem, that is to say, it is necessary to restore a solution from the scattering data. This method, also known as the nonlinear Fourier transform (NFT), has recently begun to attract much attention in areas where NLSE is used to describe various types of optical signals. In particular, NFT found use in telecommunication applications, in which it offered a new method for compensating for the effects acting on a signal during its propagation in optical communication lines~\cite{Yousefi2014III, Le2014, Gui2017a, Wahls2017, Gui2018, Civelli2019}. The method is also used to describe and analyze various physical phenomena~\cite{Gelash2019PRL, Mullyadzhanov2019}.
In this regard, interest in numerical methods for solving the direct and inverse problem has increased. Reviews of existing methods are given in~\cite{Yousefi2014II, Turitsyn2017Optica, Vasylchenkova2019a}.
Separately, it is worth mentioning the fast methods for calculating NFT, which have less asymptotic complexity compared to traditional approaches, and are called the fast nonlinear Fourier transform, FNFT~\cite{Wahls2013, Wahls2015, Vaibhav2018, Wahls2018, Chimmalgi2019}.

The solution to the direct problem consists of finding a discrete and continuous spectrum of a differential operator of a special form. If the signal is specified by $M$ values, then the continuous spectrum can be effectively determined using FNFT in $\mathcal{O}(M\log^2 M)$ operations, however, finding a discrete spectrum remains a difficult computational task requiring at least $\mathcal{O}(M^2)$ arithmetic operations. The FNFT-approach proposes to search for the discrete spectrum as the eigenvalues of a special matrix, which requires $\mathcal{O}(M^2)$ operations, although, for a large number of points in the signal, the discrete spectrum can be found for $\mathcal{O}(M\log^2 M)$ operations by down-sampling the signal without significant loss of accuracy~\cite{Chimmalgi2019}. The disadvantages of the method include the need to work with large matrices, as well as the need to filter the found eigenvalues.
The TIB (Toeplitz inner bordering)~\cite{Frumin2015_TIB} method, which has the complexity $\mathcal{O}(M^2)$, also belongs to the matrix methods for finding the discrete spectrum.

Other common approach for determining a discrete spectrum is to use iterative methods for finding zeros of the Jost function~$a(\zeta)$, which are the discrete eigenvalues of considered operator. In problems with a large discrete spectrum, the random search may be more efficient than the matrix methods~\cite{Yousefi2014II}, but if the initial approximations are unsuccessful, the speed of the method can slow down significantly. One of the variants of the search algorithm was proposed in~\cite{Burtsev1998} and consisted in the fact that the search area was divided into a grid with a large step, at the nodes of which the value of the desired function was calculated. The size of the discrete spectrum must first be determined through the argument variation of the scattering data calculated from the continuous spectrum. If it was not possible to determine all discrete eigenvalues on a large grid, then one needs to select a grid with a smaller step and to repeat the process. It should be noted here that determining the size of the discrete spectrum by the argument variation can be difficult when the change in the value of the argument occurs over a sufficiently large interval that exceeds the region of the continuous spectrum that is interesting from a practical point of view.

Another important class of algorithms for finding the discrete spectrum are methods based on the calculation of contour integrals in the complex plane. This approach was first proposed in~\cite{DelvesLyness1967}, and it was applied to the direct scattering problem in~\cite{Vasylchenkova2018a}. The idea of the method is to use the Cauchy theorem and compose, by using Newton's identities a polynomial whose roots coincide with the roots of the function~$a(\zeta)$ inside a given contour. Since the construction of such a polynomial in the case of a large number of roots is unstable, it is proposed to divide the search region into subdomains. In such a case, a situation is possible when the contour is near the zero of the function~$a(\zeta)$, which leads to an increase of the integration error.
After determining the approximate values of the discrete spectrum, it is proposed to refine them using the iterative method, for example, Newton's method.
Generally speaking, most of the discrete spectrum search algorithms discussed give values that can be refined using iterative methods; therefore, this approach is preferably applied whenever possible.
The contour integral method works well when the size of the discrete spectrum is small, and the discrete spectrum itself can be well localized initially (for example, in the case of data transmission via optical lines, when a small number of pulse types are operated on).

Among the algorithms for zero search of analytic functions, one can note the bisection-exclusion method~\cite{Yakoubsohn2005}, in which the search region is divided into square regions, and regions that do not contain function zeros are discarded using the proposed criterion.

Finally, we note the method based on minimizing the functional constructed using the continuous spectrum, and also taking into account the knowledge about the size of the discrete spectrum~\cite{Aref2019}. The method works fast with a small size of the discrete spectrum but may require the selection of a good initial approximation. With a large number of discrete eigenvalues, the method often converges to local minima.

In this paper we propose a new method for finding discrete eigenvalues for the direct Zakharov-Shabat problem based on complex plane along special trajectories leading to eigenvalues.
It is proposed to use jumps of the argument of the function~$a(\zeta)$ for such trajectories, the localization of which does not require great accuracy.
This method is less prone to computational errors. For its operation, it is not necessary to calculate the continuous spectrum, but if the continuous spectrum is calculated, then the information from it can be used to speed up the method. Also, the method does not require a priori knowledge of the size of the discrete spectrum. In the case of the correct selection of parameters for setting up the method, it allows finding all discrete eigenvalues taking into account their multiplicity faster than matrix methods (including FNFT) and contour integrals. The motivation for developing this method arose when it was discovered that the problem of finding a large-sized discrete spectrum for complex pulses (tens and hundreds of discrete eigenvalues) could not be correctly solved by existing algorithms in an acceptable time. This algorithm can find a discrete spectrum several orders of magnitude faster.

\section{The direct Zakharov-Shabat problem}

The propagation of an optical signal in a fiber under certain conditions is described by the NLSE, which can be written in dimensionless units as
\begin{equation}
 i \frac{\partial q}{\partial z} + \frac{1}{2} \frac{\partial^2 q}{\partial t^2} + \sigma |q|^2 q = 0 {,}
 \label{eq:nlse}
\end{equation}
where $q(t,z)$ is a complex-valued function that describes an optical signal depending on the time coordinate $t$ and spatial $z$.

The main idea of the method proposed by Zakharov and Shabat in 1972~\cite{ZakharovShabat1972} was the transition from a nonlinear equation to an operator equation for which Eq.~(\ref{eq:nlse}) is a compatibility condition.
This method is called the inverse scattering problem method (ISM), also known as the nonlinear Fourier transform (NFT).
The direct problem of Zakharov-Shabat consists in solving the spectral problem and can be written out as a system
\begin{equation}
 \left\{
 \begin{aligned}
 -\partial_{t} \psi_1 + q(t,0) \psi_2 = i\zeta \psi_1 \\
 \partial_{t} \psi_2 -\sigma q^{*}(t,0) \psi_1 = i\zeta \psi_2 \\
 \end{aligned}
 \right. {,}
 \label{eq:zs}
\end{equation}
where $\psi_{1,2}$ are vector components of the eigenfunction, $q(t,0)$ is initial ($z=0$) optical signal distribution, $\zeta$ is the spectral parameter. The coefficient $\sigma = \pm 1$ corresponds to focusing and defocusing cases. A set of $\{\zeta_k\}$ satisfying the problem~(\ref{eq:zs}) for a given distribution of $q(t, 0)$ is called a ``nonlinear spectrum'' consisting of a continuous and discrete parts.
A continuous spectrum exists for any initial parameters and is located on the entire real axis, while the existence of a discrete spectrum depends on the optical field distribution, and the eigenvalues lie in the upper complex half-plane.
Physically the discrete spectrum corresponds to the solitons existing in the initial signal.

In addition to the nonlinear spectrum, when solving the problem~(\ref{eq:zs}), it is necessary to determine the so-called ``scattering data'' --- the coefficients $a(\zeta)$ and $b(\zeta)$,
which together with the spectrum completely determine the initial potential $q(t,0)$. The coefficients are determined by the formulas
\begin{eqnarray}
 a(\zeta) = \lim_{t \to +\infty} \psi_1(t,\zeta) e^{i\zeta t} {,}\nonumber \\
 b(\zeta) = \lim_{t \to +\infty} \psi_2(t,\zeta) e^{-i\zeta t}
\end{eqnarray}
For a continuous spectrum located on the whole real line,
the reflection coefficient is determined by the formula:
\begin{equation}
 \rho(\zeta) = \frac{b(\zeta)}{a(\zeta)}{,} \quad \zeta \in \mathbb{R} {.}
\end{equation}
Zeros of the coefficient $a(\zeta_k)$ in the upper complex half-plane determine the discrete spectrum $\{\zeta_k\}$, $k = 1 ... K$, where $K$ is the number of discrete eigenvalues, and the parameter
\begin{equation}
 \rho_k = \frac{b(\zeta_k)}{a'(\zeta_k)} {,} \quad \text{where} \quad
 a'(\zeta_k) = \frac{\partial a(\zeta)}{\partial \zeta}|_{\zeta = \zeta_k} {,}
 \label{eq:scat_sol}
\end{equation}
is a phase coefficient.

As mentioned above, the nonlinear spectrum and scattering coefficients completely determine the original signal. In this case, if we know the nonlinear spectrum and scattering data, then we know the nonlinear dynamics of the signal propagating in an optical fiber.

\section{Method description}

The proposed method of finding the discrete spectrum is based on the argument principle of meromorphic function~$f$, saying that the change in the argument $\arg f(z)$ in traversing along some closed contour without self-intersections $\partial G$, bounding a simply connected region~$G$ and not passing through the zeros and poles of the function~$f$ is expressed in terms of the number of zeros~$K$ and the poles~$P$ of this function as follows:
\begin{equation}
\label{argPrinciple}
 Z - P = \frac{1}{2 \pi i} \int\limits_{\partial G} \frac{df}{f} = \frac{1}{2 \pi} \Delta_{\partial G} \arg f.
\end{equation}
Since the scattering coefficient~$a(\zeta)$ has no poles in the upper complex half-plane ($a(\zeta)$ is an analytic function), Eq.~\eqref{argPrinciple} is transformed into a simple form:
\begin{equation}
\label{argPrincipleA}
 \Delta_{\partial G} \arg f = 2 \pi Z.
\end{equation}
Thus, when choosing a suitable region~$G$ containing all $Z$ zeros of the coefficient~$a(\zeta)$, the change in the argument turns out to be $2\pi Z$. It follows that, on the contour~$\partial G$, the main value of the argument~$\arg a(\zeta)$, taking values within $(-\pi; \pi]$, must have at least $K$ discontinuities of the first kind. Due to the analyticity of the coefficient~$a(\zeta)$, there are $K$ jumps in its argument that continuously extend into the domain along the trajectories~$\{\gamma_k \}$ and disappear at the points corresponding to zeros~$\{ \zeta_k \}$ of the coefficient~$a(\zeta)$, in which its argument is not defined.

As a result, we propose to find numerically the jump points of the argument~$\arg a(\zeta)$ on the contour that spans all zeros of~$\{\zeta_k \} $, and then track the jumps to the corresponding zeros. For a complete solution to the problem of finding the discrete spectrum of the ZSP, the following sequence of steps must be performed:
\begin{enumerate}

 \setlength\itemsep{0.1em}

 \item Determining the boundaries of the domain containing a discrete spectrum.

 \item Identification of jumps in the argument of the coefficient $a(\zeta)$ at the boundary of the domain containing a discrete spectrum.

 \item Choosing a step to track the trajectories of the argument jump and refinement of the jump points of the argument at the domain boundary.

 \item Tracking phase jump trajectories.

 \item Refinement of discrete eigenvalues using iteration methods.

\end{enumerate}

We now discuss each step in more detail.

\subsection{Determining the boundaries of the domain containing a discrete spectrum}

If it is necessary to find only the discrete spectrum, then the region~$[L_{\xi}, R_{\xi}]$ on the real axis, including all the real parts of the discrete spectrum, is selected from the linear Fourier transform of the signal~$q(t)$. We propose using the following criterion: if $\tilde{q} (\omega) $ is the linear Fourier transform of the original signal, and $\tilde{q}_{\max} = \max\limits_{\omega} |\tilde{q}(\omega) |^2$, then $L_{\xi}$ is chosen as the minimum frequency at which the following relation is satisfied:
\begin{equation}
 |\tilde{q}(L_{\xi})|^2 / \tilde{q}_{\max} = C_q,
\end{equation}
where the constant $C_q$ was chosen as $10^{-4}$.
The value $R_{\xi}$ is chosen in the same way as for the maximum frequency at which this relation is satisfied. Such choice of the region limited the real parts of the discrete spectrum well in all the considered examples.


The upper bound~$U$ of the region of searching for zeros of the coefficient~$a(\zeta)$ is determined based on the Parseval equation for the energy of the initial signal
\begin{equation}
\label{ParsevalRelation}
 E_t = \int\limits_{-\infty}^{\infty} \left| q(t) \right|^2 dt = 4 \sum\limits_{k=1}^K \eta_k - \frac{1}{\pi} \int\limits_{-\infty}^{\infty} \log \left| a(\xi) \right|^2 d\xi = E_d + E_c.
\end{equation}
and is chosen as
\begin{equation}
\label{UpperBound}
 U = \min[1.1 \cdot 0.25 (E_t - E_c), \ \log(0.9 L_{\text{MAX}}) / T]
\end{equation}
where $L_{\text{MAX}}$ is the maximum number represented by the type of real number used on the computer. If the continuous spectrum has not been calculated before, then the energy~$E_c$ is omitted in the estimate~\eqref{UpperBound}. Resulting search area is $G = [L_{\xi}, R_{\xi}] \times [0, U]$.

\subsection{Identification of jumps in the argument of the coefficient $a(\zeta)$ at the boundary of the domain containing a discrete spectrum}

The determination of the points~$ \{g_k \} $ at which the argument~$ \arg a(\zeta) $ jumps does not require a large number of calculations of the coefficient~$ a(\zeta) $ at the boundary~$ \partial G $ and can be done on a fairly coarse grid. The only thing that needs to be noted is that the step on the section of the border passing through the real axis ($ [L_{\xi}, R_{\xi}] $) must still be smaller than the step on other boundaries, which is automatically performed in the case of preliminary finding the continuous spectrum. The algorithm for determining the jumps of the argument used by us is presented below (algorithm~\ref{alg:jumpDetection}).

\begin{algorithm}[H]
\label{alg:jumpDetection}
\SetAlgoLined
\textbf{Input:} $\{z_j\} \in \partial G$ --- ordered nodes on the boundary, $\{a(z_j)\}$, $j=0...J$.

\textbf{Output:} $\{g_k\}$ --- phase jumps.\newline

\For{j = 1...J}{

 If{ $\arg a(z_j) \cdot \arg a(z_{j-1}) < 0$ \textbf{and} $\left| \arg a(z_j) - \arg a(z_{j-1}) \right| > 1.3 \pi$
 }{
 \text{add} $(z_{j} + z_{j - 1}) / 2$ \text{to} $\{g_k\}$
 }
}{}{}
 \caption{Phase jump detection on the boundary $\partial G$}
\end{algorithm}

After identifying the points~$ \{g_k \} $ on the segment $ [L_{\xi}, R_{\xi}] $, it is proposed to choose the step to bypass the rest of the boundary~$\partial G$ equal to the minimum distance between these points:
\begin{equation}
\label{boundStep}
 h_{\partial G \backslash [L_{\xi}, R_{\xi}]} = \min\left(0.5 \min\limits_{k,j}(|g_k - g_j|), 0.01 (U + R_{\xi} - L_{\xi})\right), \ \{g_k\} \in [L_{\xi}, R_{\xi}].
\end{equation}

\subsection{Choosing a step to track the trajectories of the argument jump}

In this paper, it is proposed to choose a step to bypass the trajectories~$\{\gamma_k\} $ as a constant, although, of course, various approaches to adaptive step selection can be applied here. In our calculations, we used the step
\begin{equation}
\label{traectStep}
 h_{\gamma} = C_h \min\limits_{k,j}(|g_k - g_j|), \ \{g_k\} \in \partial G,
\end{equation}
after determining which, for all jumps of the argument~$ \{g_k \} $ found earlier on the boundary~$ \partial G $, using the bisection method, they were refined up to $ 0.5 h_{\gamma} $. The choice of the constant $C_h$ depends on the type of discrete spectrum. In all further calculations, unless otherwise specified, the value $C_h = 1/15$ was used.

\subsection{Tracking phase jump trajectories}

Argument trajectories can be tracked in various ways. Here, we give for demonstration purposes only one of the options, perhaps not the most effective. The tracking method in some aspects is similar to the well-known method for isolating contours in images called ``marching squares''~\cite{Lorensen1987, Maple2003}.

It is proposed to bypass the trajectories, finding at each step a pair of points in the domain~$G$ located on different sides of the trajectory. We will only consider the start of the algorithm from breakpoints located on the real axis. Start from other segments of the boundary is carried out similarly with the difference in orientation on the plane.
To determine the coefficient~$a(\zeta)$ value, we propose to use the fastest known method of the second order of accuracy for the ZS system --- the Ablovitz-Ladik (AL) method.

\begin{algorithm}[htbp]
\label{alg:tracking}
\SetAlgoLined
\textbf{Input:} Phase jumps $\{g_k\}$, initial jump orientations $\{\varphi_k\}$ ($\varphi_k = 0$ is lower boundary, $\varphi_k = \pi/2$ is right boundary, $\varphi_k = \pi$ is upper boundary, $\varphi_k = 3\pi/2$ is left boundary), $k = 1...K$, and the step $h_{\gamma}$.

\textbf{Output:} \{$\zeta_k$\}.\newline

\For{k = 1...K}{
 $l_k^0 = g_k - 0.5 h_{\gamma} \exp\left(i \varphi_k \right) $\;
 $r_k^0 = g_k + 0.5 h_{\gamma} \exp\left(i \varphi_k \right) $\;


 \While{true}{
 $\text{shift} := h_{\gamma} \exp\left(i (\varphi_k + \pi / 2)\right)$\;
 $ r_k^* := r_k^m + \text{shift} $\;
 \eIf{$\arg a(r_k^{m}) \cdot \arg a(r_k^*) < 0$ \textbf{and} $\left| \arg a(r_k^{m}) - \arg a(r_k^*) \right| > 1.3 \pi$}{
 $ l_k^{m+1} := r_k^*$\;
 $ r_k^{m+1} := r_k^{m} $\;
 $\varphi_k := \varphi_k - \pi/2$\;
 }{
 $ l_k^* := l_k^{m} + \text{shift}$\;

 \eIf{$\arg a(l_k^*) \cdot \arg a(r_k^*) < 0$ \textbf{and} $\left| \arg a(l_k^*) - \arg a(r_k^*) \right| > 1.3 \pi$}{
 $ l_k^{m+1} := l_k^*$\;
 $ r_k^{m+1} := r_k^*$\;
 }{
 \eIf{$\arg a(l_k^{m}) \cdot \arg a(l_k^*) < 0$ \textbf{and} $\left| \arg a(l_k^{m}) - \arg a(l_k^*) \right| > 1.3 \pi$}{
 $ l_k^{m+1} := l_k^{m}$\;
 $ r_k^{m+1} := l_k^* $\;
 $\varphi_k := \varphi_k + \pi/2$\;
  }{
  $\zeta_k := 0.25 \cdot \left( l_k^{m} + r_k^{m} + l_k^* + r_k^* \right)$\;
  \textbf{break}\;
  }
 }
 }
 }
}{}{}
 \caption{Tracking phase jumps from the boundary $\partial G$}
\end{algorithm}

The points $\{ l_k^0 := g_k - 0.5 h_{\gamma},\ r_k^0 := g_k + 0.5 h_{\gamma} \}$ are selected as starting. One step along the path~$\gamma_k$ is to define a new pair of points between which the desired trajectory passes.
To do this, one need first find the value of the argument~$\arg a(\zeta)$ at the point~$r_k^*:= r_k^0 + h_{\gamma}$. The criterion for jump determination:
\begin{equation}
\label{breakCriterium}
 \arg a(r_k^{0}) \cdot \arg a(r_k^*) < 0, \quad \left| \arg a(r_k^{0}) - \arg a(r_k^*) \right| > 1.3 \pi
\end{equation}
If it is met, the new pair of points $\{ l_k^1 := r_k^*,\ r_k^1 := r_k^0 \}$ is selected.
If not, then the argument~$\arg a(\zeta)$ is calculated at the point~$ l_k^* := l_k^0 + h_{\gamma} $.
If the trajectory places between the points $l_k^*$ and $r_k^*$, then they are selected as the next pair: $ l_k^1 := l_k^*$ and $ r_k^1 := r_k^*$.
If the trajectory is between the points $l_k^0$ and $l_k^*$, then we set $\{l_k^1 := l_k^0, \ r_k^1 := l_k^* \}$. If, however, it is not possible to determine the presence of a jump between all three pairs using the criterion~\eqref{breakCriterium}, then we assume that the jump reached the point of the discrete spectrum $\zeta_k$, therefore we stop the process of moving along the trajectory. The full form of the algorithm is presented below (algorithm~\ref{alg:tracking}).


\subsection{Refinement of discrete eigenvalues}

Since at the last step the discrete eigenvalue was determined only up to a square with side~$h_{\gamma}$, then it was necessary to refine it. We performed it using the Muller method~\cite{Muller1956}
\begin{eqnarray}
\label{Muller}
 &&x_{n+1} = x_n - 2 f_n / d, \\ \nonumber
 &&d = \max \left[ w - \sqrt{w^2 - 4 f_n g}, \ w + \sqrt{w^2 - 4 f_n g} \right], \\ \nonumber
 &&w = f[x_{n-1}, x_n] + f[x_{n-2}, x_n] - f[x_{n-2}, x_{n-1}], \\ \nonumber
 &&g = f[x_{n-2}, x_{n-1}, x_n], \nonumber
\end{eqnarray}
having the order of convergence $\approx 1.84$ (which was higher than the order of the secant method), and having a sufficiently large area of convergence~\cite{Vasylchenkova2019a}. In particular, the method can converge from real starting points. In our case $f = a(\zeta)$, and we solve the equation $a(\zeta) = 0$. The Muller method uses divided differences:
\begin{equation}
 f[x_1, x_2] = \frac{f(x_2) - f(x_1)}{x_2 - x_1}, \ f[x_k, ..., x_{k+m}] = \frac{f(x_{k+1}, ..., x_{k+m}) - f(x_k, ..., x_{k+m-1})}{x_{k+m} - x_k}.
\end{equation}
In our case the middle of the square and 2 points close to it were chosen as the starting points for the Muller method. The use of the Muller method is also due to its speed: at each iteration, it is only necessary to calculate the value of the function $a(\zeta)$ once, unlike, for example, the Newton method, where the derivative must also be calculated. On the other hand, it may be preferable to use the Newton method on coarse grids. It should be noted that at this stage it was necessary to abandon the AL method since its approximation order did not allow reliable determination of the values of the discrete spectrum. Instead, we used the previously proposed exponential ES4 scheme of the 4th order approximation~\cite{Medvedev2020_OE}, based on the of Magnus decomposition~\cite{Magnus1954}. Although it should be noted here that the operating time of all the considered 4th-order schemes is tens of times longer than the AL method operation time, therefore, the step with refining the discrete spectrum can take a significant part of the entire method operation time.

\subsection{Some remarks}

Note that a case is possible when the argument jump trajectories enter and exit the search region of the discrete spectrum without passing through the zeros of the coefficient $a(\zeta)$. This can happen, for example, when one of the trajectories leading to a zero enters the region several times. A large number of such trajectories can slow down the algorithm.
However, the proposed algorithm does not need to start from all jump points of the argument $\{g_k\}$, but only from those for which the derivative of the argument is positive to the left and the right of the jump point in the positive direction of the contour $ \partial G $ tracing. This is due to the fact that the negative derivative of the argument in the vicinity of the jump indicates that this jump belongs to a trajectory that has already entered the domain $G$ earlier.

Note that the method does not involve the use of any a priori information about the number of zeros of the coefficient $a (\zeta) $. One can make sure that all zeros were found by the method by checking Parseval equality~\eqref{ParsevalRelation}, as well as the remaining conservation integrals~\cite{Yousefi2014II}.

We also note that there is an approach where first all points of the argument jump are determined on the real axis, and then only these points extend to the upper half-plane, i.e., instead of a rectangular domain, which contains zeros of the coefficient $a(\zeta)$, a sufficiently large segment on the real axis is considered. But this approach has three significant drawbacks. Firstly, the jumps in the argument can be located much further than the domain in which the continuous spectrum is noticeably different from zero, which will be shown later on the example of an oversoliton with a large discrete spectrum. Secondly, these trajectories can go beyond the boundaries of the domain in which it is possible to calculate the value of the coefficient $a(\zeta)$ using the given type of real numbers. Finally, the trajectories from the extreme points of the jump to the zeros of the coefficient $a(\zeta)$ may be of great length, what can slow down the operation of the algorithm. The latter drawback can be corrected by an adaptive selection of the tracking step.

Other special cases, including the intersection of trajectories and the case of multiple roots, will be discussed in detail in the next section.

\section{Numerical experiments}

We now turn to the consideration of specific examples of the proposed algorithm. In all tests the time domain of the optical pulse was chosen as $[-T; T]$, which is divided by $M$ equal intervals. So the discretization step was $\tau=2T/M$, and grid nodes were $t_n=-T+\tau n$, where $n=0,...,M$. Here we calculated the continuous spectrum to obtain a more precise estimation for the upper boundary.

\subsection{Examples with oversolitons}

\begin{figure}[t]
 \includegraphics[width=1.0\linewidth]{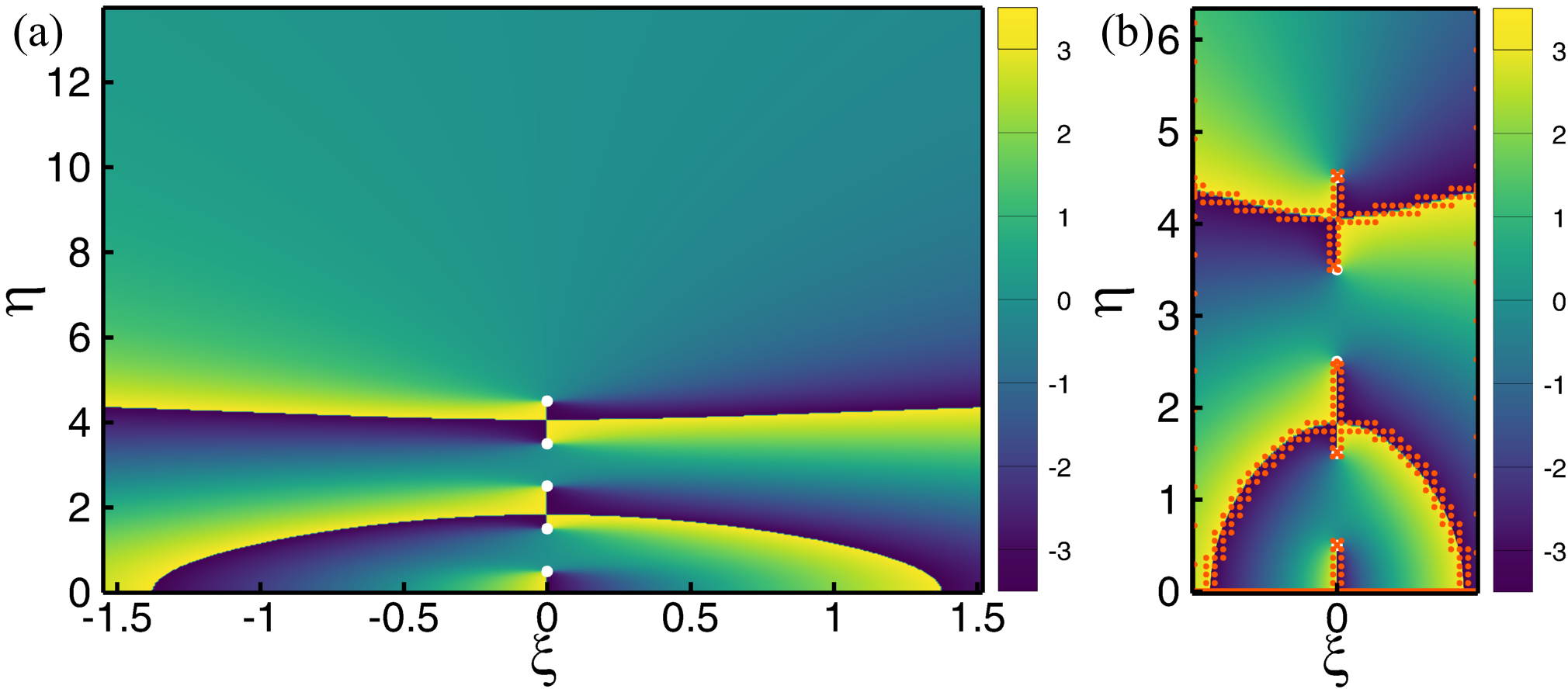}
 \vspace{-0.1cm}
 \caption{Oversoliton with parameters $A=5$, $C=0$: (a) the domain $G$, $\arg a(\zeta)$ for it and the discrete spectrum of the signal under investigation denoted by white circles; (b) enlarged central area, which depicts the operation of the algorithm (orange dots). Used $M = 2^{14}$ points of the signal.}
 \label{FIG:O_5_0_ALL}
\end{figure}

Firstly, we consider the results of applying the algorithm to chirped oversolitons
\begin{equation}
 \label{chirpedOversoliton}
 q(t) = A[\mbox{sech}(t)]^{1+iC},
\end{equation}
where $A$ and $C$ are real parameters.
The discrete spectrum $\{\zeta_k\}$, $k=\overline{0,K-1}$ of given pulses is presented as follows:
\begin{equation}\label{exact_dzeta}
\zeta_k = i\bigl(\sqrt{A^2 - C^2/4} - 1/2 - k\bigr),\quad k = 0,\ldots,[\sqrt{A^2 - C^2/4} - 1/2],
\end{equation}
where square brackets denote the integer part of the expression. Here we set $T = 30$.

Let us first consider the oversoliton with parameters $A=5$, $C=0$ (Fig.~\ref{FIG:O_5_0_ALL}) whose discrete spectrum consists of 5 eigenvalues located on the imaginary axis. The algorithm perfectly defines 5 argument jumps. As one can see, the jump trajectories in this case intersect and diverge at right angles. Since the algorithm~\ref{alg:tracking} first checks for a gap at the right point, in this case there is no collision and the algorithm moves along each path to the right side. This ensures that all zeros of the function~$a(\zeta)$ are achieved.

\begin{figure}[ht]
 \centering
 \includegraphics[width=0.3\linewidth]{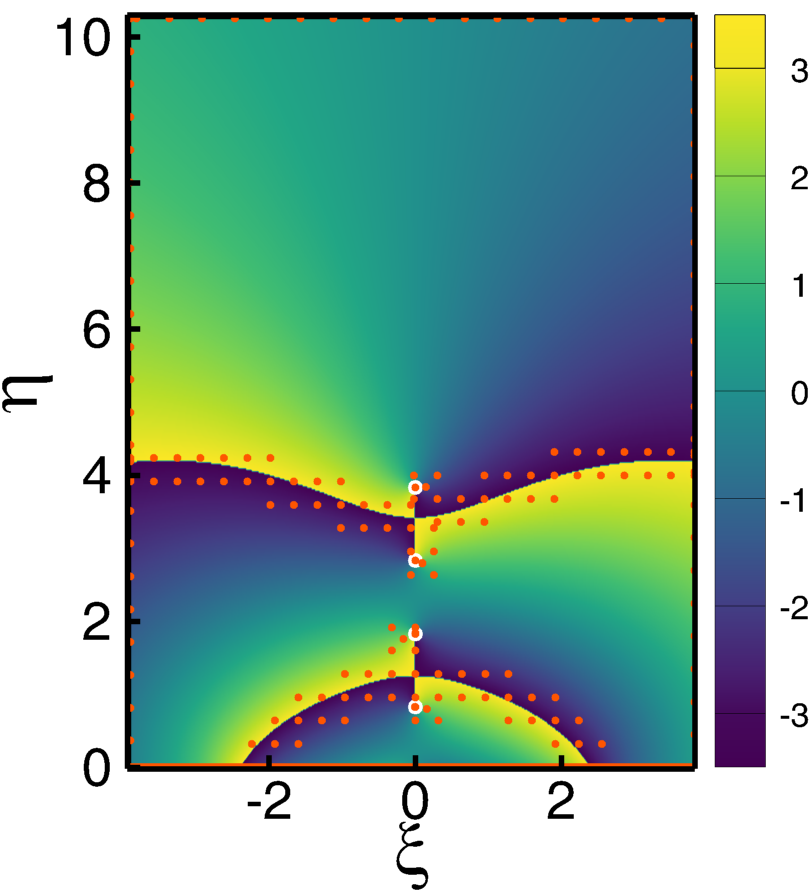}
 \vspace{-0.1cm}
 \caption{Oversoliton with parameters $A=5$, $C=5$: the domain $G$, $\arg a(\zeta)$ for it and the discrete spectrum of the signal under investigation denoted by white circles. The operation of the algorithm is also shown.}
 \label{FIG:O_5_5_ALL}
\end{figure}

If the chirp $C = 5$ is added to the given oversoliton (Fig.~\ref{FIG:O_5_5_ALL}), one eigenvalue disappears. The proposed algorithm also determines the entire discrete spectrum.

\begin{figure}[b]
 \includegraphics[width=1.0\linewidth]{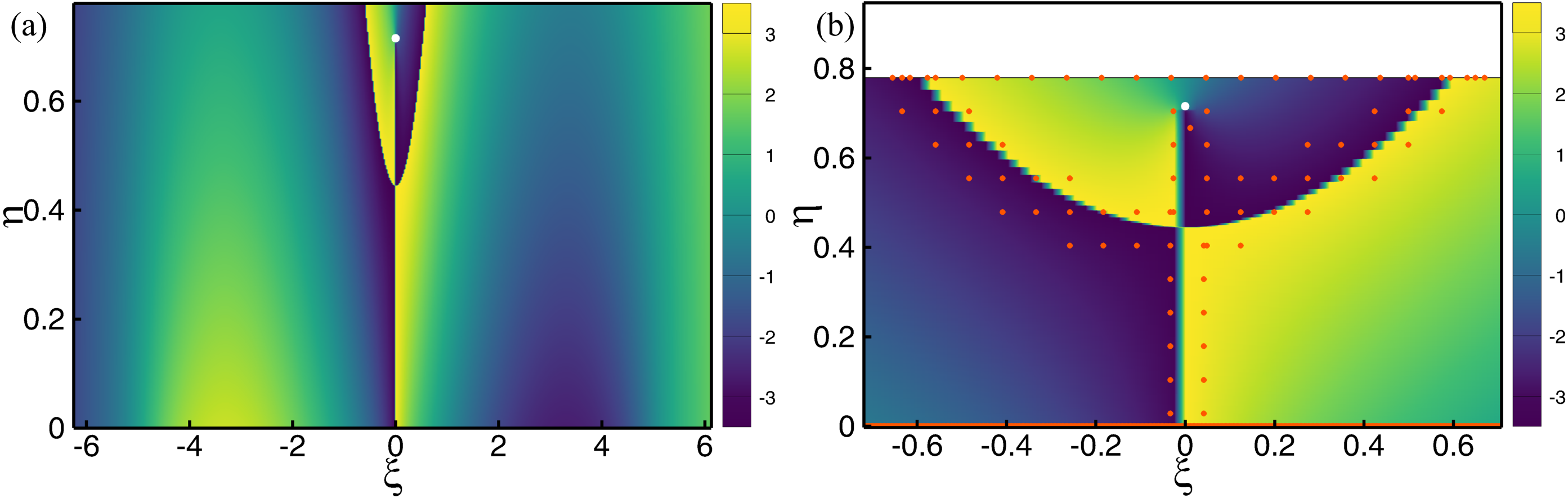}
 \vspace{-0.1cm}
 \caption{Oversoliton with parameters $A=5$, $C=9.7$: (a) the domain $G$, $\arg a(\zeta)$ for it and the discrete spectrum of the signal under investigation denoted by white circles; (b) enlarged central area, which depicts the operation of the algorithm (orange dots).}
 \label{FIG:O_5_9,7_ALL}
\end{figure}

Fig.~\ref{FIG:O_5_9,7_ALL} shows the operation of the algorithm for oversoliton with parameters $A = 5$, $C = 9.7$, which has only one discrete eigenvalue. In this case, there are 3 argument jumps at the boundary of the domain $G$, however, the algorithm does not start from the jump on the real axis, since the derivative is negative there. Thus, only one jump leads to the discrete eigenvalue, and the second to a jump on the real axis.

\begin{figure}[ht]
 \includegraphics[width=1.0\linewidth]{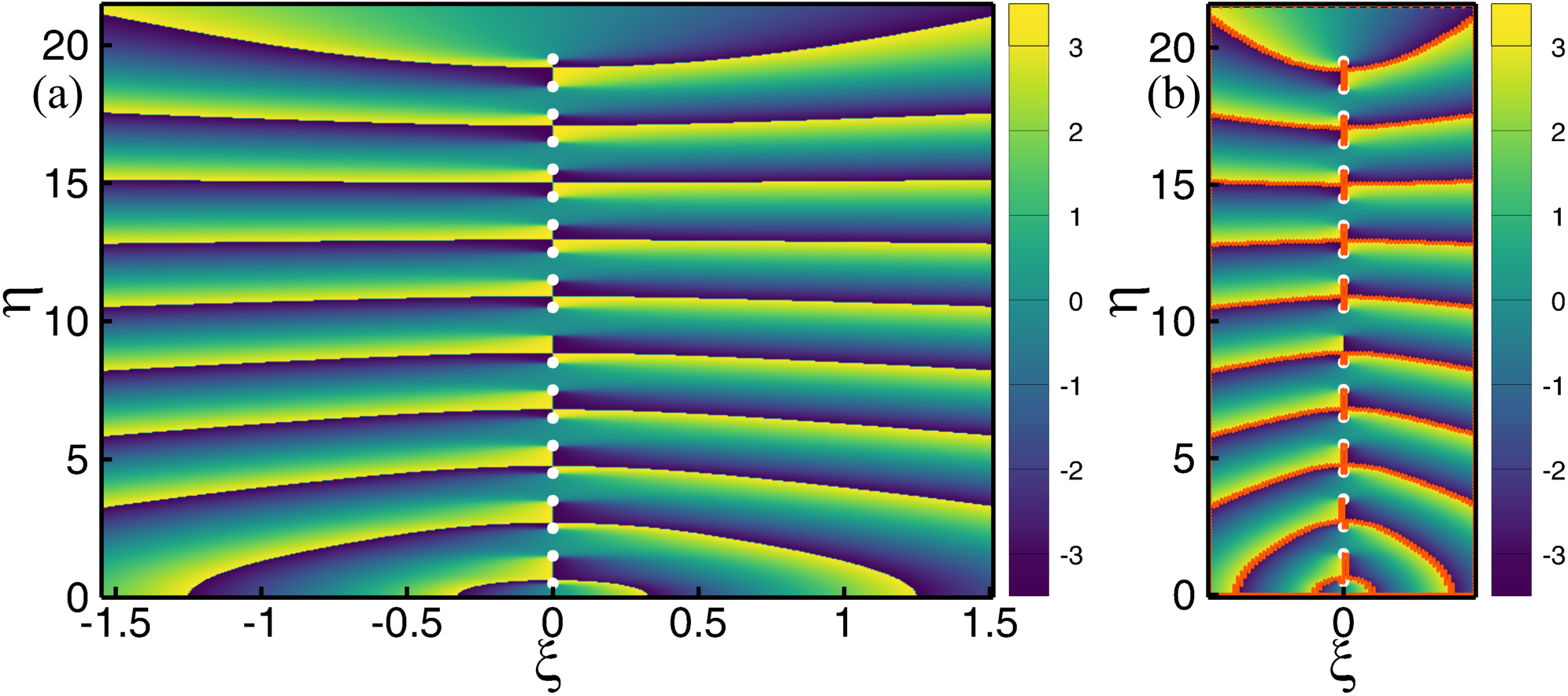}
 \vspace{-0.1cm}
 \caption{Oversoliton with parameters $A=20$, $C=0$: (a) the domain $G$, $\arg a(\zeta)$ for it and the discrete spectrum of the signal under investigation denoted by white circles; (b) enlarged central area, which depicts the operation of the algorithm (orange dots).}
 \label{FIG:O_20_0_ALL}
\end{figure}

In conclusion, we consider an oversoliton with parameters $A = 20$, $C = 0$ (Fig.~\ref{FIG:O_20_0_ALL}), which has 20 eigenvalues, the maximum value of which is close to the boundary at which the coefficient $a(\zeta)$ can still be calculated using double precision. In this case, it is more convenient to determine the jumps of the argument on a rectangular contour bounding the discrete spectrum, and not on the real axis, since the jumps appear on it far beyond the region where the continuous spectrum is nonzero. In particular, due to the restriction to a rectangular region, the length of the trajectories leading to discrete eigenvalues is reduced, which decreases the algorithm operation time.

\subsection{32-soliton solution example}

\begin{figure}[ht]
 \includegraphics[width=1.0\linewidth]{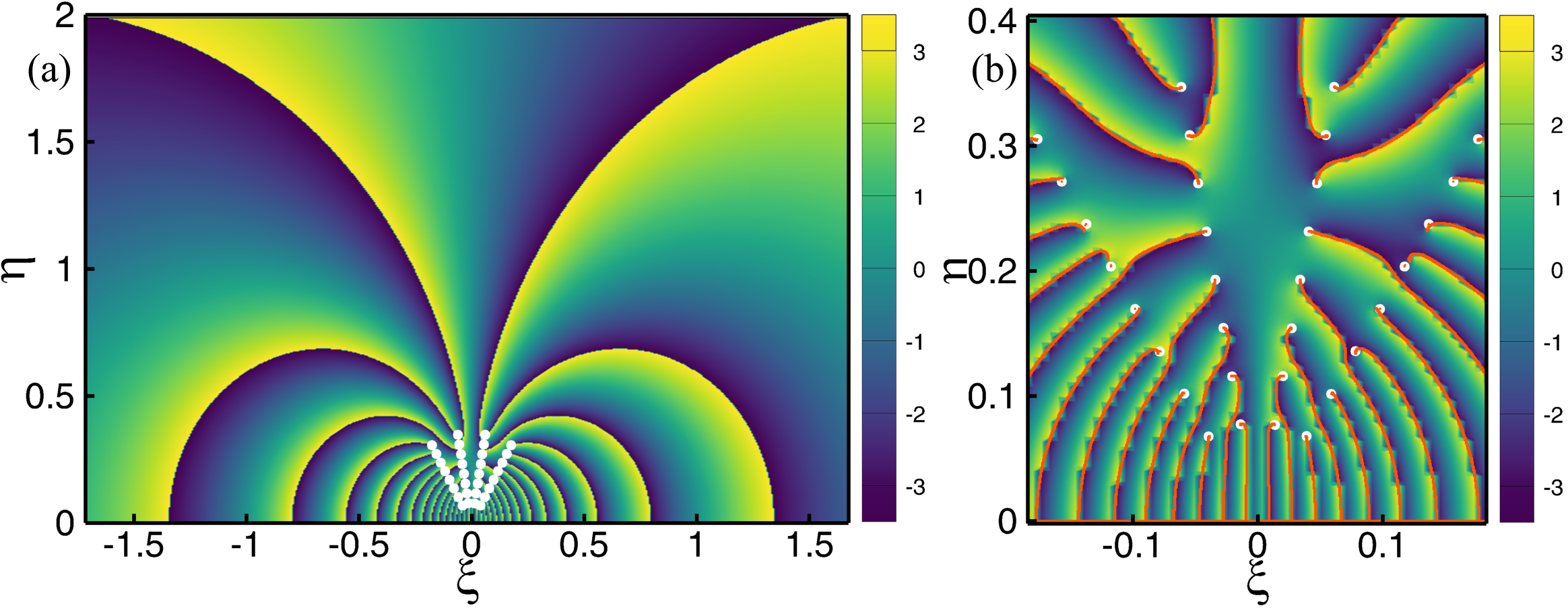}
 \vspace{-0.1cm}
 \caption{32-soliton solution: (a) the domain $G$, $\arg a(\zeta)$ for it and the discrete spectrum of the signal under investigation denoted by white circles; (b) enlarged central area, which depicts the operation of the algorithm.}
 \label{FIG:V_ALL}
\end{figure}

This example is taken from the paper~\cite{Vaibhav2017} and is a multisoliton solution constructed from the given discrete spectrum by the Darboux method~\cite{Aref2016}. The solution plot and its discrete spectrum, consisting of 32 eigenvalues, are presented in Fig.~\ref{FIG:V_ALL}. In this test, the FNFT algorithm based on the search for zeros of the polynomial specifying $a(\zeta)$ is not able to find all 32 discrete eigenvalues using any discretization $M$ of the signal. Nevertheless, the proposed algorithm determines the entire discrete spectrum much faster.

\subsection{Example with multiple eigenvalues}

\begin{figure}[ht]
 \includegraphics[width=1.0\linewidth]{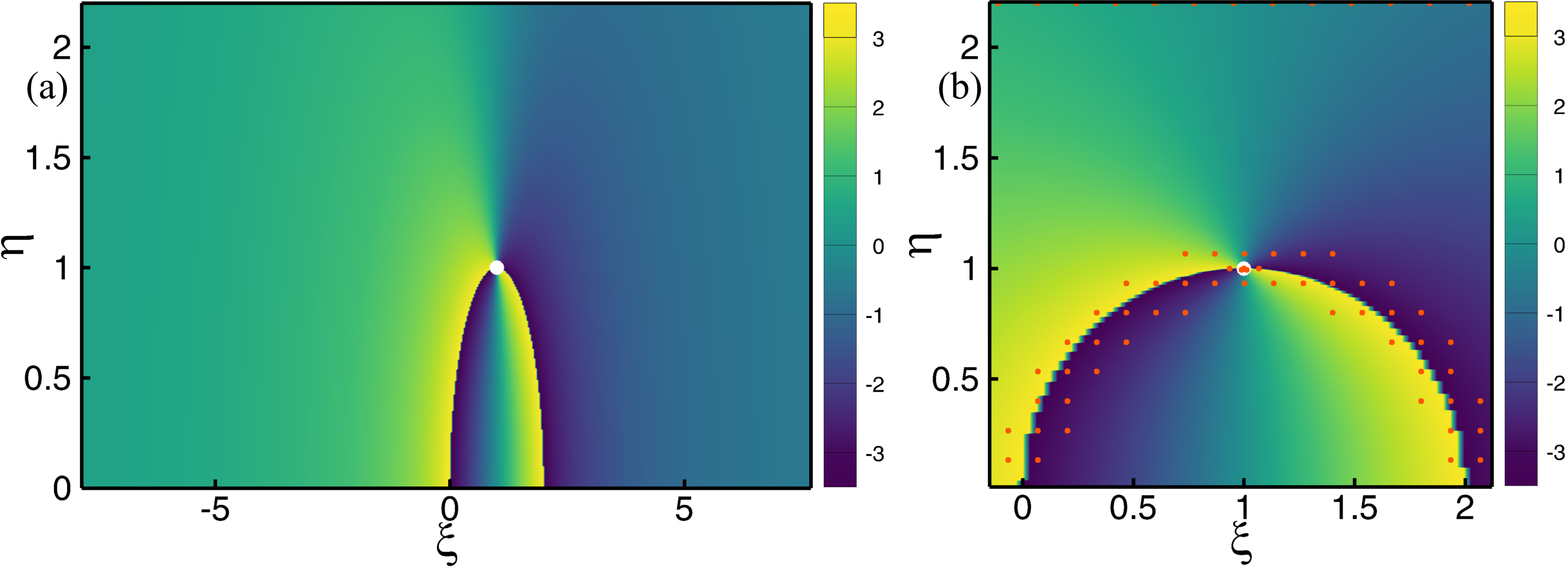}
 \vspace{-0.1cm}
 \caption{Soliton solution with one multiple discrete eigenvalue: (a) the domain $G$, $\arg a(\zeta)$ for it and the discrete spectrum of the signal under investigation denoted by white circles; (b) enlarged central area, which depicts the operation of the algorithm.}
 \label{FIG:Double_ALL}
\end{figure}

The example is taken from~\cite{Garcia-Gomez2018} and is a soliton with a discrete eigenvalue $\zeta = \xi + i\eta$ of multiplicity 2, as well as two normalization constants $Q_{11}$ and $Q_{10}$:
\begin{equation}
 q(t)=\frac{h(t)}{f(t)},
\end{equation}
where
\begin{eqnarray}
 &&h(t)=-i 4 \eta e^{-i \arg Q_{11}} e^{-i 2 \xi t} \left\{e^{-X} \left[-\left|Q_{11}\right|^{2}(2 \eta t+2)-\eta Q_{11}^{*} Q_{10}\right] \right. \\ \nonumber &&\qquad \left.+e^{X}\left[\left|Q_{11}\right|^{2}(2 \eta t)+\eta Q_{11} Q_{10}^{*}\right]\right\}, \\ \nonumber
 &&f(t)=\left|Q_{11}\right|^{2}[\cosh(2 X)+1] + 2\left|Q_{10} \eta+Q_{11}(2 \eta t+1)\right|^{2}, \\ \nonumber
 &&X= 2 \eta t-\log \frac{\left|Q_{11}\right|}{4 \eta^{2}}. \nonumber
\end{eqnarray}
In the example shown in Fig.~\ref{FIG:Double_ALL} we used $\zeta = 1 + i$, $Q_{11} = Q_{10} = 1$. As one can see, the number of phase jumps at the boundary of the domain $G$ is equal to two. There is a trajectory to a single discrete eigenvalue from each of them. Thus, the algorithm allows finding a discrete spectrum taking into account the multiplicity of eigenvalues: if $N$ trajectories have come to the same discrete eigenvalue, so its multiplicity is equal to $N$.

\subsection{Rectangular pulse example}

\begin{figure}[t]
 \includegraphics[width=1.0\linewidth]{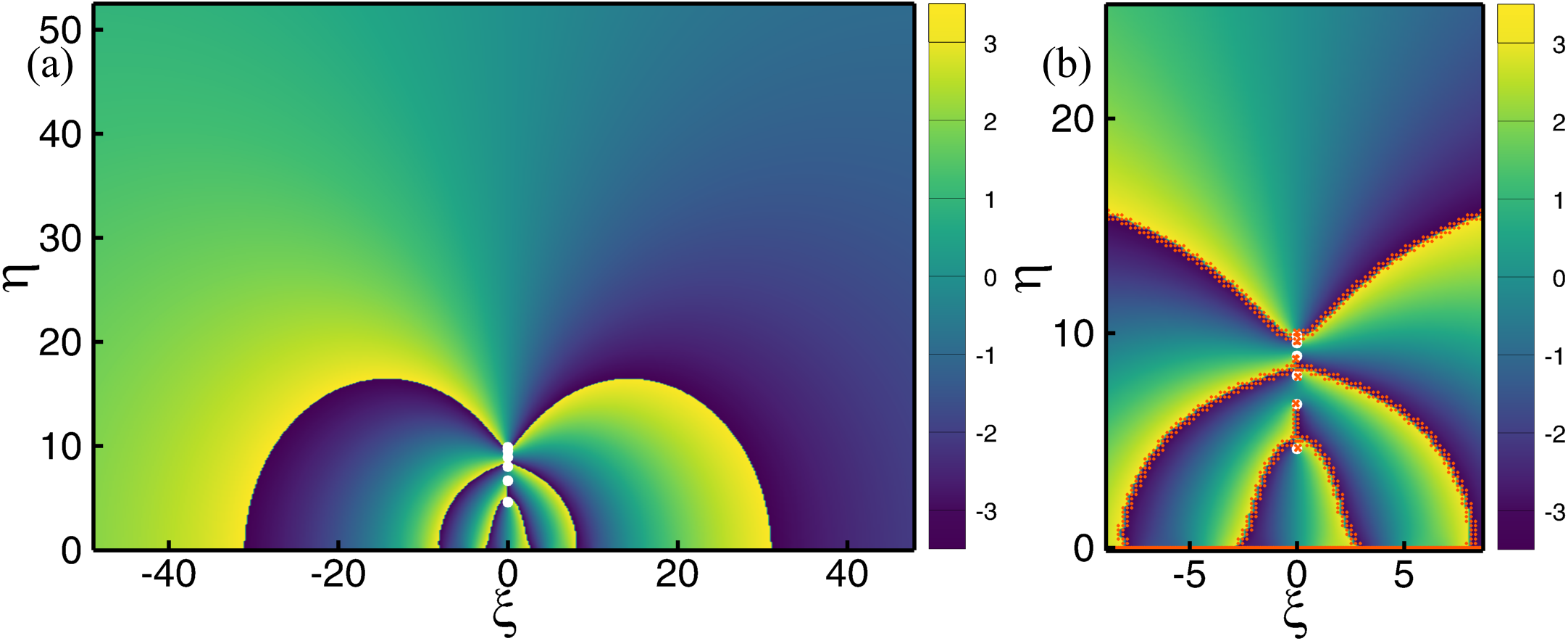}
 \vspace{-0.1cm}
 \caption{Rectangular pulse with parameters $A=10$, $T=1$: (a) the domain $G$, $\arg a(\zeta)$ for it and the discrete spectrum of the signal under investigation denoted by white circles; (b) enlarged central area, which depicts the operation of the algorithm.}
 \label{FIG:Rect_ALL}
\end{figure}

Fig.~\ref{FIG:Rect_ALL} presents the results of the algorithm operation for the rectangular pulse
\begin{equation}
 q(t)=
 \begin{cases}
 A, \ t \in [-T; T]\\
 0, \text{otherwise},
 \end{cases}
\end{equation}
the discrete spectrum size of which is equal to 6 for $A=10$, $T=1$. The eigenvalues of the given signal are concentrated to one point, therefore, the distance between adjacent discrete eigenvalues is significantly different. In this case, the previous step choice by the formula~\eqref{traectStep} is incorrect and the two upper trajectories cannot distinguish the eigenvalues with the maximum imaginary part. The trajectories come to different eigenvalues when choosing the tracking step $h_{\gamma}$ with the constant $C_h = 1/25$ at least.


\subsection{Algorithm runtime analysis and comparison with other methods}

\begin{figure}[b]
 \includegraphics[width=1.0\linewidth]{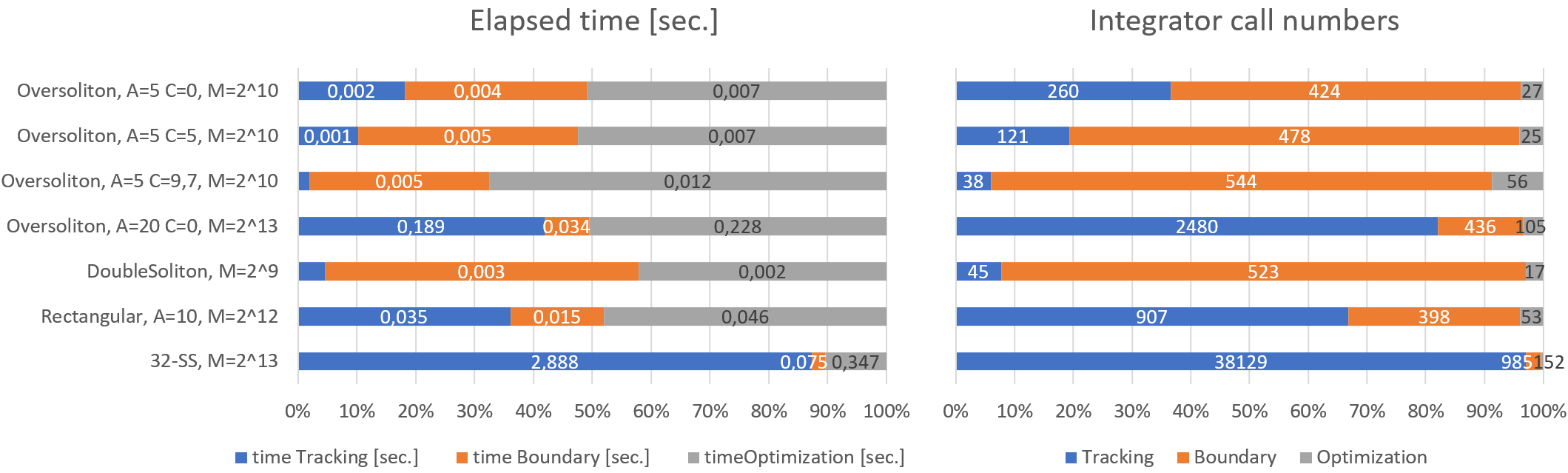}
 \vspace{-0.1cm}
 \caption{The ratio of the execution time and the number of integrator calls to solve the direct ZSP between the various steps of the method --- moving along the trajectories, calculating the function $a (\zeta) $ at the boundary of the domain $G$ and refining the discrete eigenvalues by the iterative method.}
 \label{FIG:AlgTimes}
\end{figure}

The main stages of the proposed method, which is spent most of the time on, are the calculation of the value of $a (\zeta) $ at the boundary of the discrete spectrum search domain, the movement along trajectories (tracking) and iterative refinement. Fig.~\ref{FIG:AlgTimes} for the previously reviewed tests shows the relationship between the execution time of the main stages of the algorithm, as well as the relationship between the number of integrator calls (AL and ES4) at different stages. For each test, we consider the minimum number of points per signal~$M$, at which all discrete eigenvalues were localized using the proposed algorithm.

In this case it is assumed that the continuous spectrum was not calculated in advance; therefore, the value~$a(\zeta)$ must also be determined on the segment of the real axis~$ [L_{\xi}, R_{\xi}] $. The step size to bypass this part of the boundary~$\partial G$ was chosen according to the criterion presented in~\cite{Medvedev2020_OE}: we took the same number of points $N_{\xi} = M$ in the spectral domain, defined a spectral step as $d\xi = \pi/(2L)$, and the size of the spectral interval as $L_{\xi} = \pi/(2\tau)$. Such choice of the step made it possible to localize all breaks of the argument. As a result, a total of 400-500 calls to the AL method were required to bypass the border for the tests considered. The tracking stage required a smaller number of calls of the AL method in tests with a discrete spectrum size of no more than 5. In the case of a 32-soliton solution, this stage occupied the bulk of the time the algorithm worked due to the large length of the trajectories (87 \% of the time and 96 \% of the total the number of integrator calls).

The step of iterative refinement in all cases, except for the test with a 32-soliton solution, took about half the total time of the algorithm. This is mainly due to the use of the ES4 scheme, as well as the fact that iterations stopped when the accuracy of $\varepsilon = 10^{-14}$ was reached. There was no limit on the number of iterations.
We note, however, that if the AL method is used instead of ES4 for the iterative refinement, the ratio of execution time between stages will be similar to the ratio between the number of calls of integrators, the refinement stage will take no more than 8\% of the elapsed time.
We also note that the Newton method in these examples works more stable in the case of small values of~$M$ than the Muller method. With its help, it was possible to localize all the discrete eigenvalues at $M$, which was 2 times less, in tests with oversolitons with parameters $ A = 5 $, $ C = 0 $ and $ A = 5 $, $ C = 9.7 $, as well as in the test with multiple eigenvalues.

We now turn to compare the algorithm with other methods of searching for discrete spectrum all calculations were performed on a single core of the Intel Core i5-9600K processor with a frequency of 4.7 GHz.

We assume that the problem of finding the discrete spectrum is solved incorrectly if at least one discrete eigenvalue was defined incorrectly. This approach differs from what was considered, for example, in~\cite{Chimmalgi2019}, where the difference between the numerically found spectrum and the exact one was gradually penalised. We adhere to this approach since one of our main goals is to find a numerical method that is guaranteed to determine the entire discrete spectrum of a signal. This is important in problems of reconstructing a field from a discrete spectrum when an incorrect determination of even one discrete eigenvalue can make a noticeable perturbation.

Our implementation of the contour integrals (CI) method was based on the second-order trapezoidal rule, as well as the AL method, which was used to calculate the value $a (\zeta)$ on the contour $\partial G$. A recursive division into four subdomains was implemented, in each of which the number of discrete eigenvalues was estimated by the argument principle. The division into subdomains continued until there were no more than four discrete eigenvalues in the new subdomain, after which approximate values of discrete eigenvalues located in this subdomain were calculated using Newton's identities. The case when the boundary dividing the region into new subdomains passed close to a discrete eigenvalue was not specially processed. Instead, in the tests, the initially left boundary of the domain $G$ was chosen asymmetrically concerning the right boundary ($ L_\xi = - 1.1 R_\xi $). After that it was checked that no one partition passes through discrete eigenvalues. With a new partition into subdomains, the values of the function $a (\zeta) $ from the previous iteration were not used, but were calculated again, although the step along the boundary of the subdomain was kept equal to the initial step for the entire contour $\partial G$. This approach turned out to be much easier to implement, however, reusing values from the previous step can give acceleration up to two times.

In testing only the discrete spectrum was calculated by all methods, therefore, for the proposed method and contour integrals, a rough estimate for the upper boundary $U$ of the discrete spectrum search domain was used.

Among the methods of the FNFT library for comparison we chose fast implementation of the Ablovitz-Ladik method and fast implementation of the Boffetta-Osborne method with splitting the exponent according to a fourth-order accuracy formula. In the first case the transfer matrix is a polynomial of the 1st degree from the exponential of the spectral parameter, and in the second it is a polynomial of the 2nd degree. Consideration of other FNFT schemes with a higher degree of a polynomial for the transfer matrix does not make much sense since the operating time of such schemes far exceeds the operating time of the PJT method in the examples considered.
Also, downsampling of the signal was applied to reduce the asymptotic complexity of the FNFT methods to $\mathcal{O}(M \log^2 M) $ although this approach could be applied to the proposed PJT method and CI method. Without downsampling FNFT method are significantly slower than PJT method (several orders of magnitude), so we did not consider this case.
For a more correct comparison of methods we did not use iterative refinement using the Newton method and the Boffetta-Osborne method, implemented in the FNFT software package itself. Instead, all methods used the same iterative refinement that we implemented based on the Newton method and the ES4 scheme. The Newton method was used instead of the Muller method, because nevertheless in these tests it worked more accurately, albeit more slowly. Accuracy of discrete eigenvalue determination was chosen equal to $ 2\cdot 10^{-14}$.

We propose to compare methods in their ability to find such initial approximations for discrete eigenvalues that various iterative methods can successfully converge. Therefore, in the method of contour integrals (CI) the minimal partition of the contour $ \partial G$ was chosen, in which this method can find adequate initial approximations for the Newton method, at which it can converge to the correct values of the discrete spectrum.

\begin{figure}[hbtp]
 \includegraphics[width=0.5\linewidth]{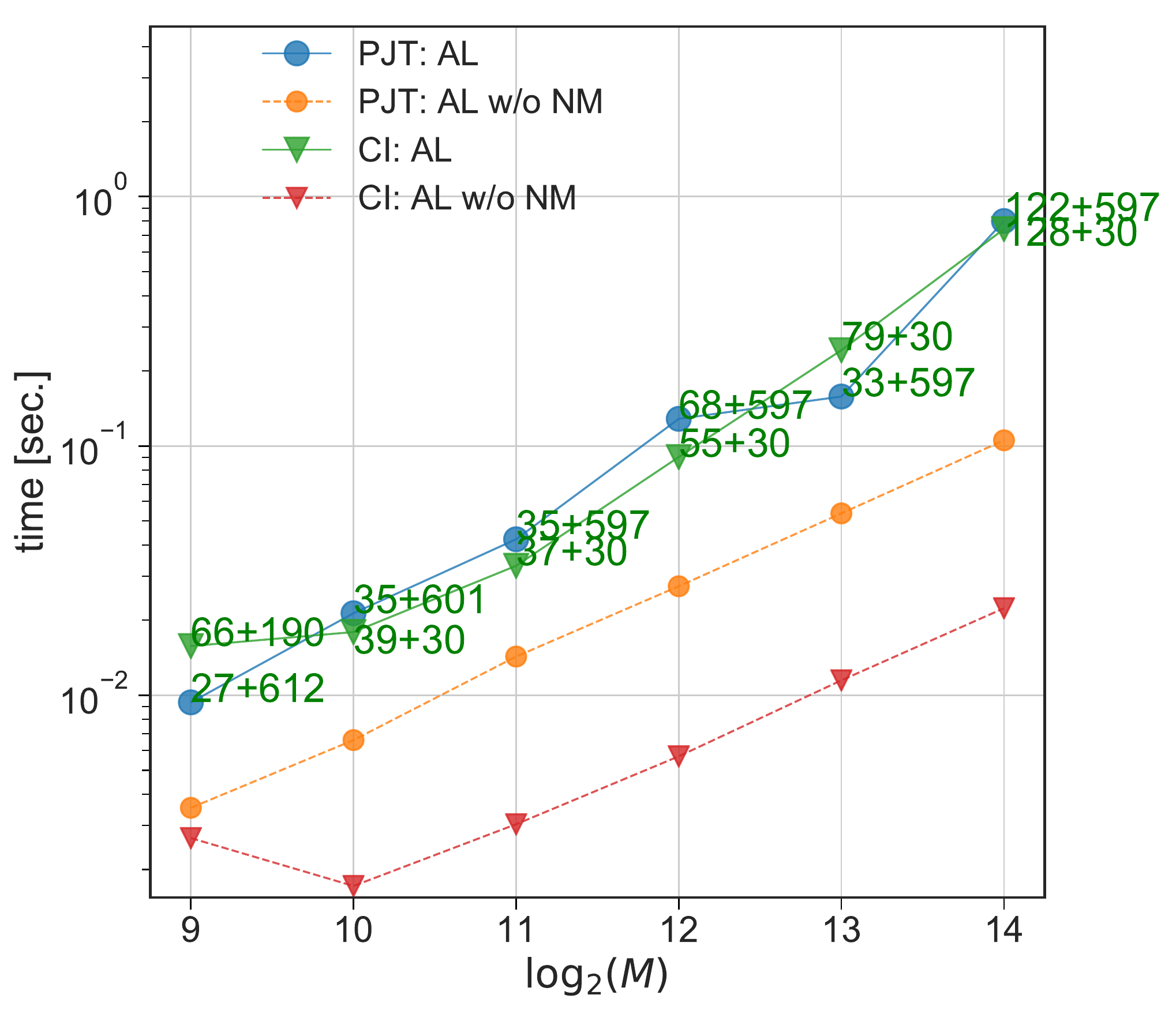}
 \includegraphics[width=0.465\linewidth]{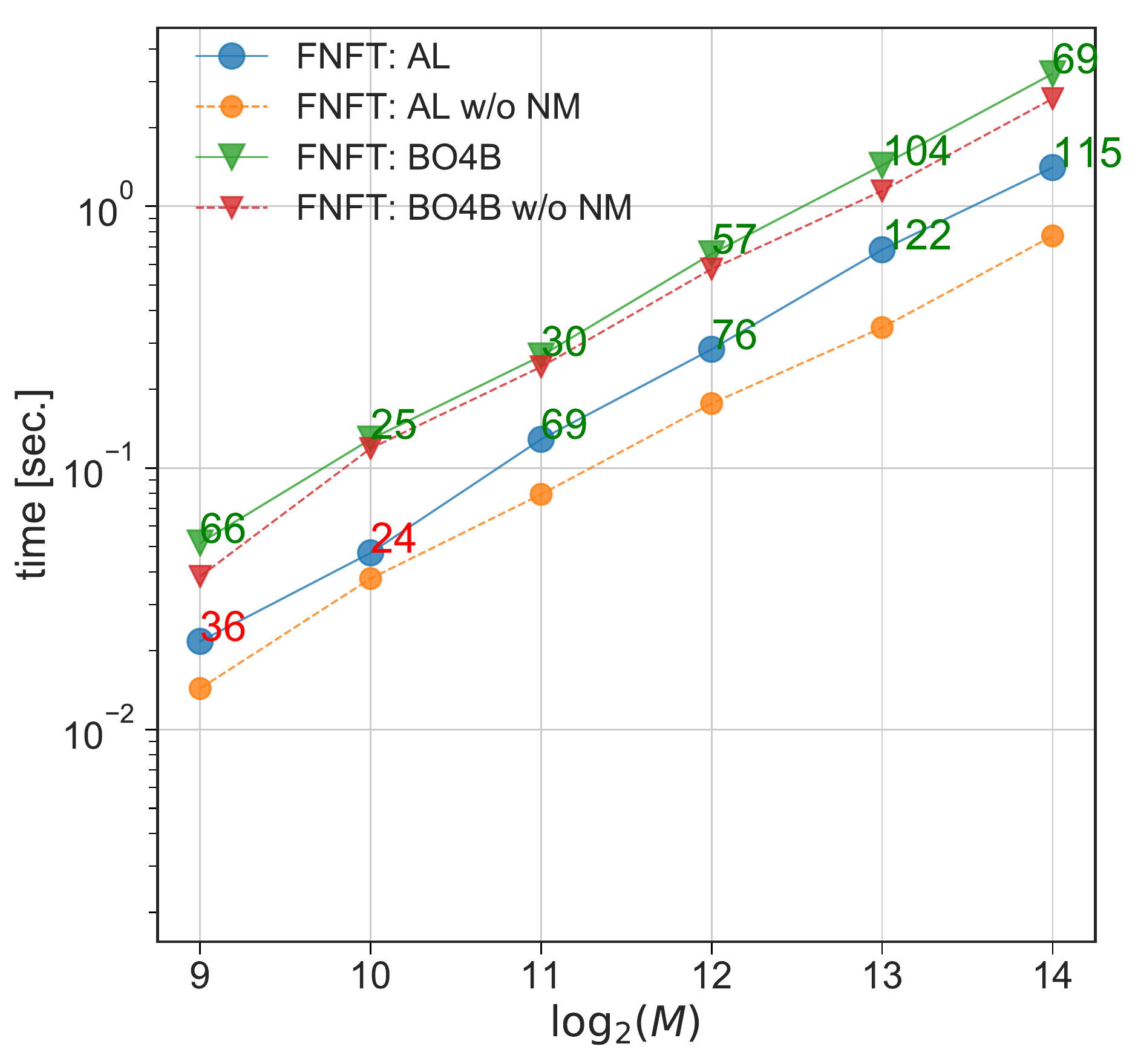}
 \vspace{-0.1cm}
 \caption{The runtime of various algorithms on the example of oversoliton with parameters $ A = 5 $ and $ C = 5 $. Solid lines indicate the total operating time, taking into account the iterative refinement of the discrete spectrum by the Newton method; dashed lines indicate the operating time without taking into account the iterative refinement. The left numbers represent the total number of calls to Newton's method. The right numbers on the plots indicate the total number of calls to the AL method for PJT and the CI method. The green colour of the numbers indicates that in this test case the discrete spectrum was correctly determined, and the red colour indicates the opposite.}
 \label{FIG:All_5_5}
\end{figure}

\begin{figure}[hbtp]
 \includegraphics[width=0.5\linewidth]{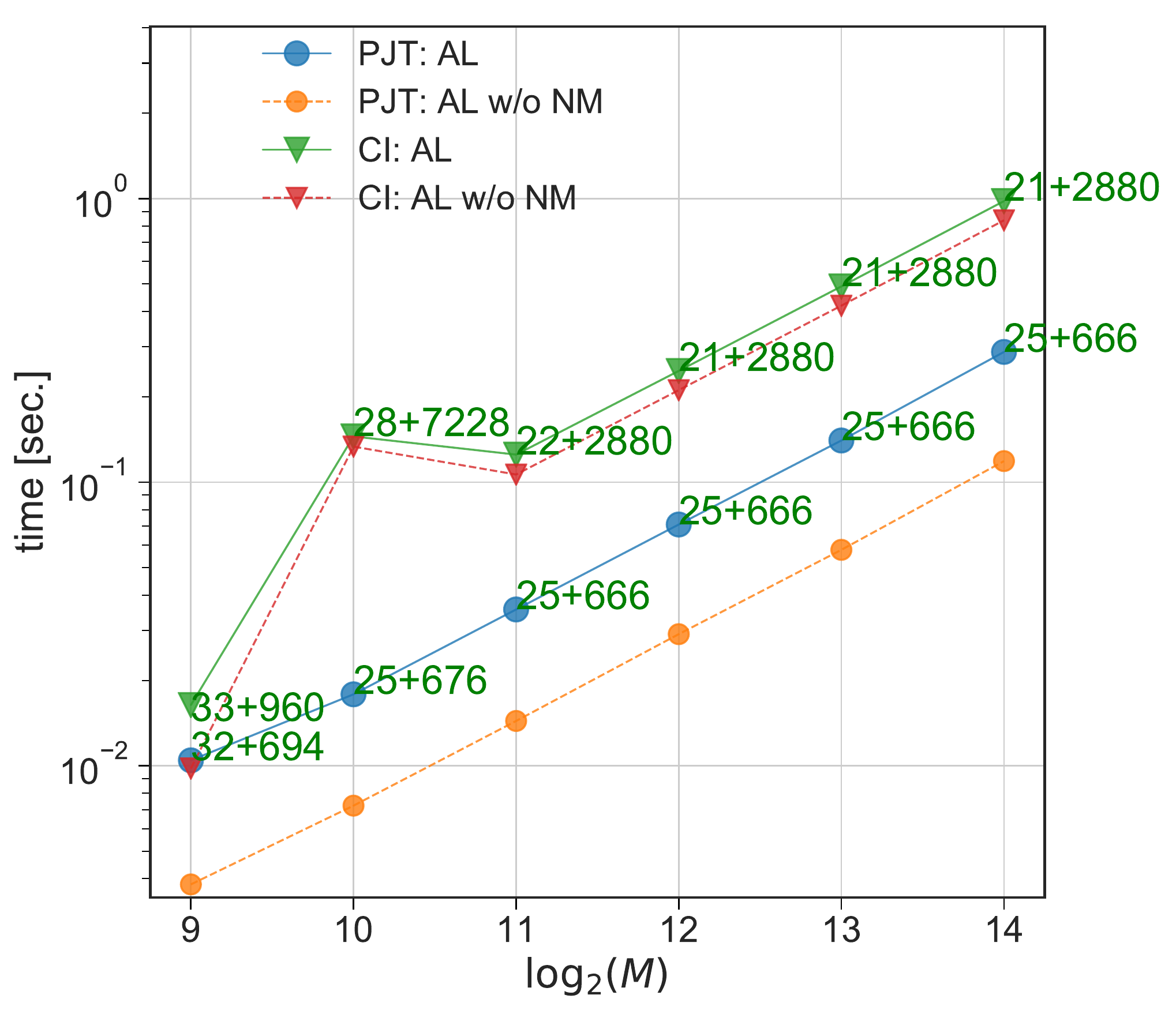}
 \includegraphics[width=0.46\linewidth]{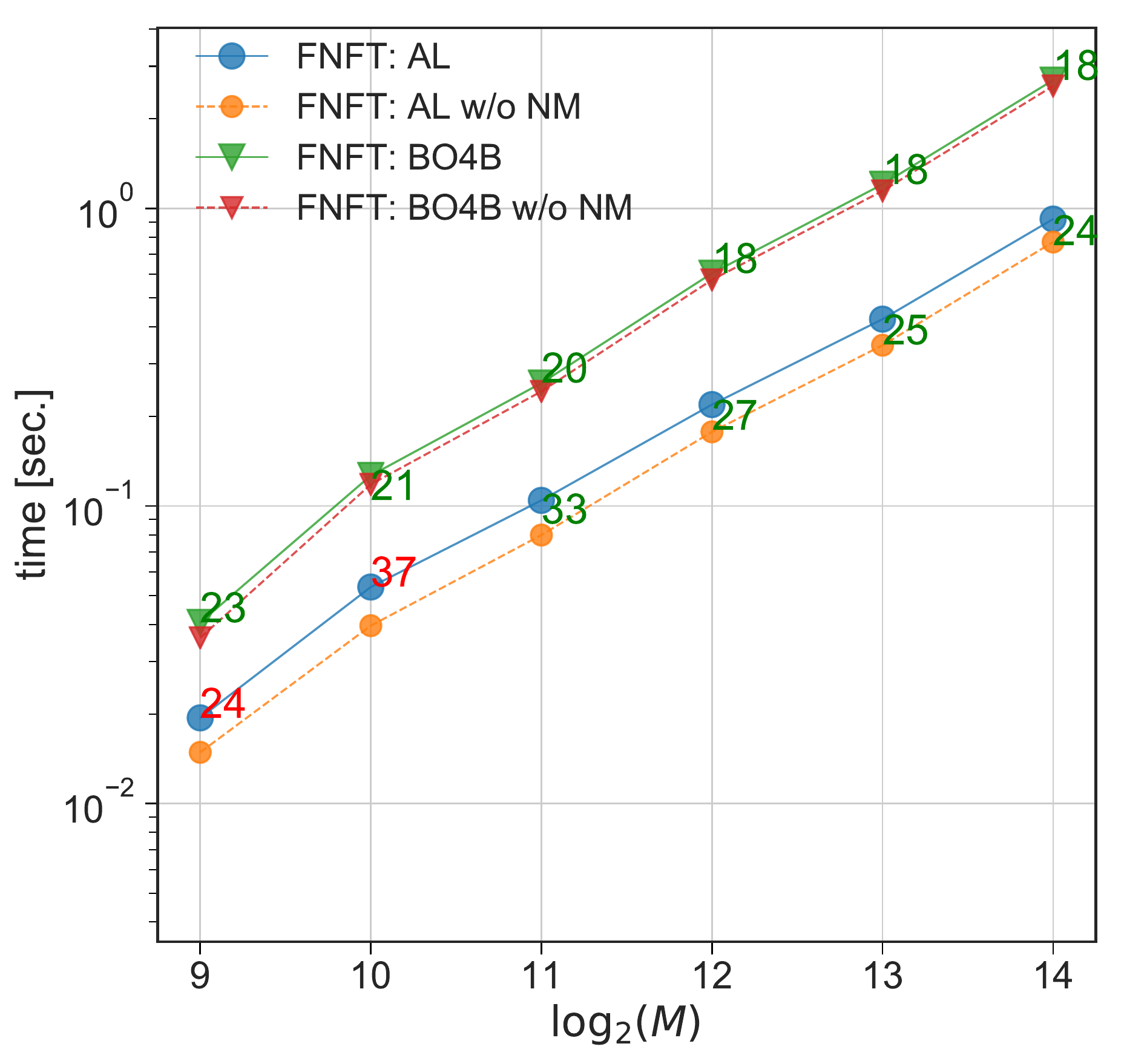}
 \vspace{-0.1cm}
 \caption{The runtime of various algorithms on the example of oversoliton with parameters $ A = 5 $ and $ C = 0 $. Solid lines indicate the total operating time, taking into account the iterative refinement of the discrete spectrum by the Newton method, dashed lines indicate the operating time without taking into account the iterative refinement. The left numbers represent the total number of calls to Newton's method. The right numbers on the plots indicate the total number of calls to the AL method for PJT and the CI method. The green colour of the numbers indicates that in this test case the discrete spectrum was correctly determined, and the red colour indicates the opposite.}
 \label{FIG:All_5_0}
\end{figure}

Figs.~\ref{FIG:All_5_5} and~\ref{FIG:All_5_0} show the results of comparing all the methods for the oversolitons with the parameters $A = 5$ and $C = 5$, as well as $ A = 5 $ and $ C = 0 $. The plots show the running time of the algorithms depending on the number of points in the signal. Also, the operating time of the methods is separately indicated without taking into account the time spent on iterative refinement by Newton's method. The discrete spectrum of oversoliton with the parameters $ A = 5 $ and $ C = 5 $ consists of 4 discrete eigenvalues, therefore, in this case, the CI method does not require a recursive partition of the contour $ \partial G $. Thus, this method is the fastest of all the considered algorithms, if you do not take into account the operating time of the Newton method. The right numbers on the plots for the CI method and the PJT method indicate the total number of individual calls to the Ablovitz-Ladik method, as well as the total number of calls to the Newton method (left numbers). As can be seen, in the case of an oversoliton with parameters $ A = 5 $ and $ C = 5 $, the CI method needs 30 calls of the AL method to determine the spectrum with sufficient accuracy so that from the found approximate values of the discrete eigenvalues Newton's method could converge. The PJT method requires about 600 points for this. The number of calls of the Newton method in both cases is approximately equal, which indicates a similar accuracy in calculating the initial approximations. Note that both FNFT methods considered are slower, and the FNFT method based on the AL scheme could not correctly determine the discrete spectrum for $ M = 2^9 $ and $ M = 2^{10} $, so the corresponding number of calls on the graph is marked in red.

If the size of the discrete spectrum is such that the CI method has to perform a recursive partition of the initial domain, as in the case of the oversoliton with parameters $ A = 5 $ and $ C = 0 $, then the number of calls to the AL method increases significantly and reaches 2880, which is more than four times the number of calls for the PJT method. Since both methods require approximately the same number of calls to the Newton method, the PJT method is three times faster than the CI method. The FNFT methods in this example also cannot compete with the PJT method in terms of speed.

\begin{figure}[hbtp]
 \includegraphics[width=0.5\linewidth]{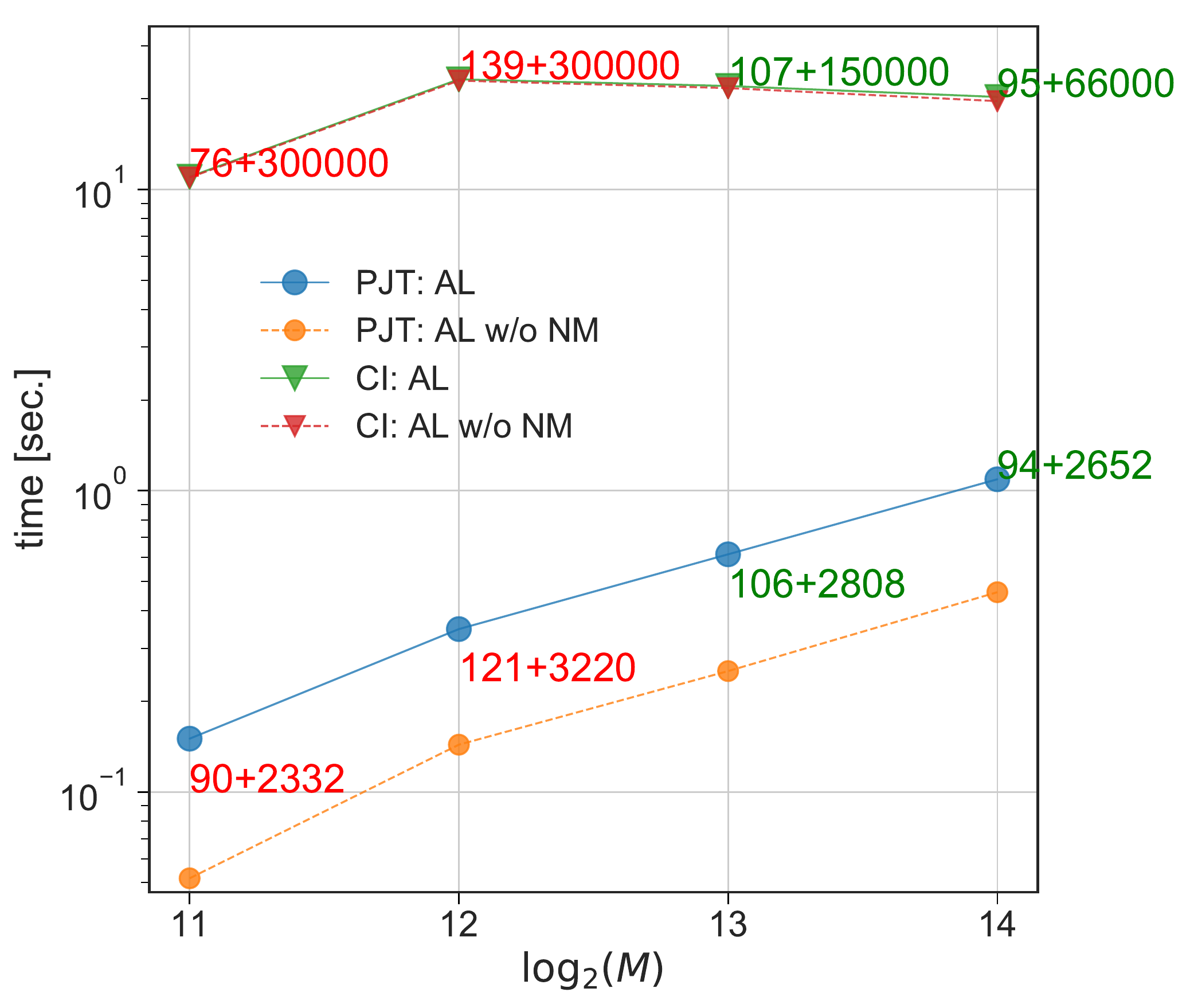}
 \includegraphics[width=0.455\linewidth]{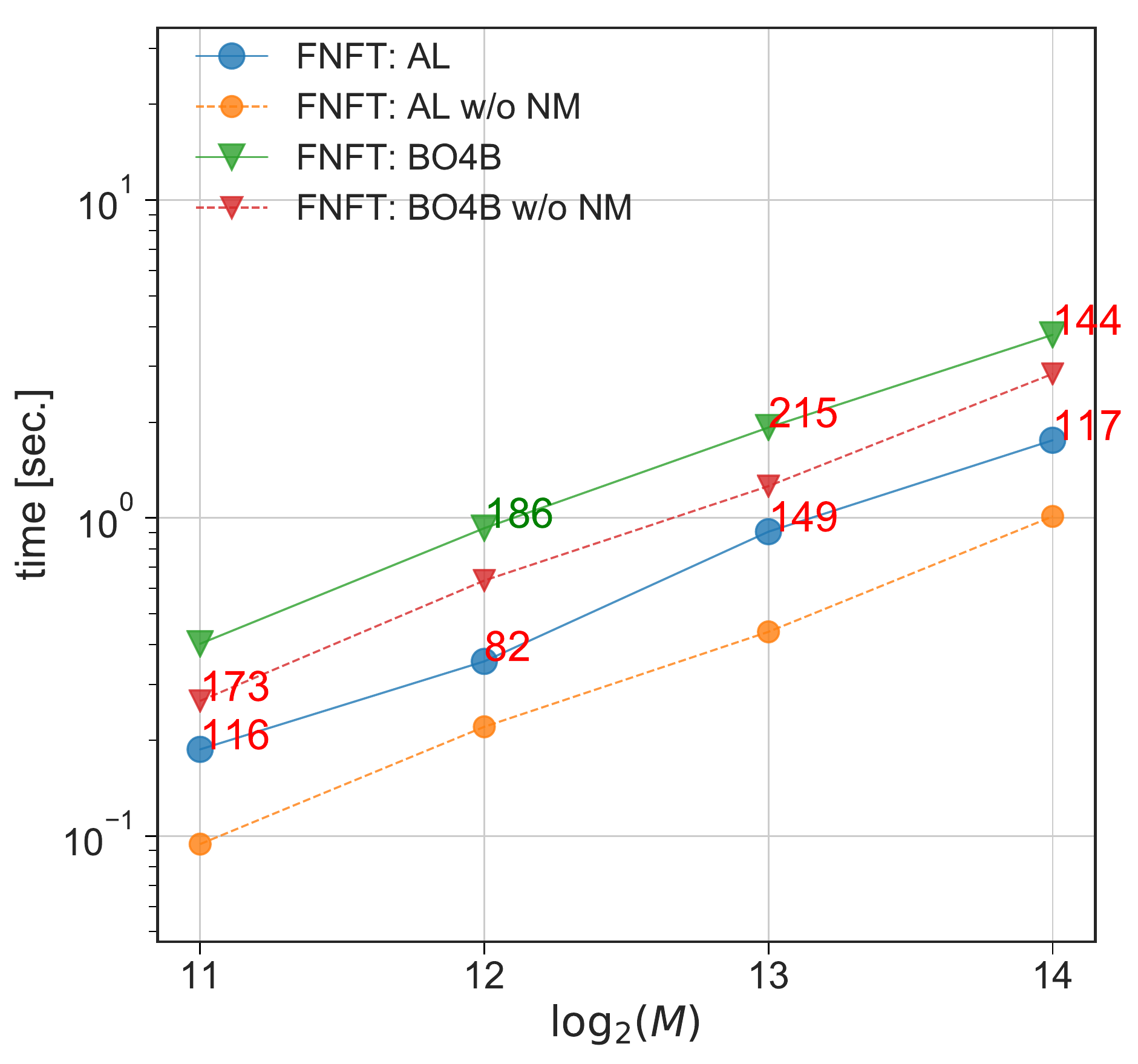}
 \vspace{-0.1cm}
 \caption{The runtime of various algorithms on the example with the oversoliton with parameters $ A = 20 $ and $ C = 0 $. Solid lines indicate the total operating time, taking into account the iterative refinement of the discrete spectrum by the Newton method, dashed lines indicate the operating time without taking into account the iterative refinement. The left numbers represent the total number of calls to Newton's method. The right numbers on the plots indicate the total number of calls to the AL method for PJT and the CI method. The green colour of the numbers indicates that in this test case the discrete spectrum was correctly determined, and the red colour indicates the opposite.}
 \label{FIG:All_20}
\end{figure}

\begin{figure}[hbtp]
 \includegraphics[width=0.5\linewidth]{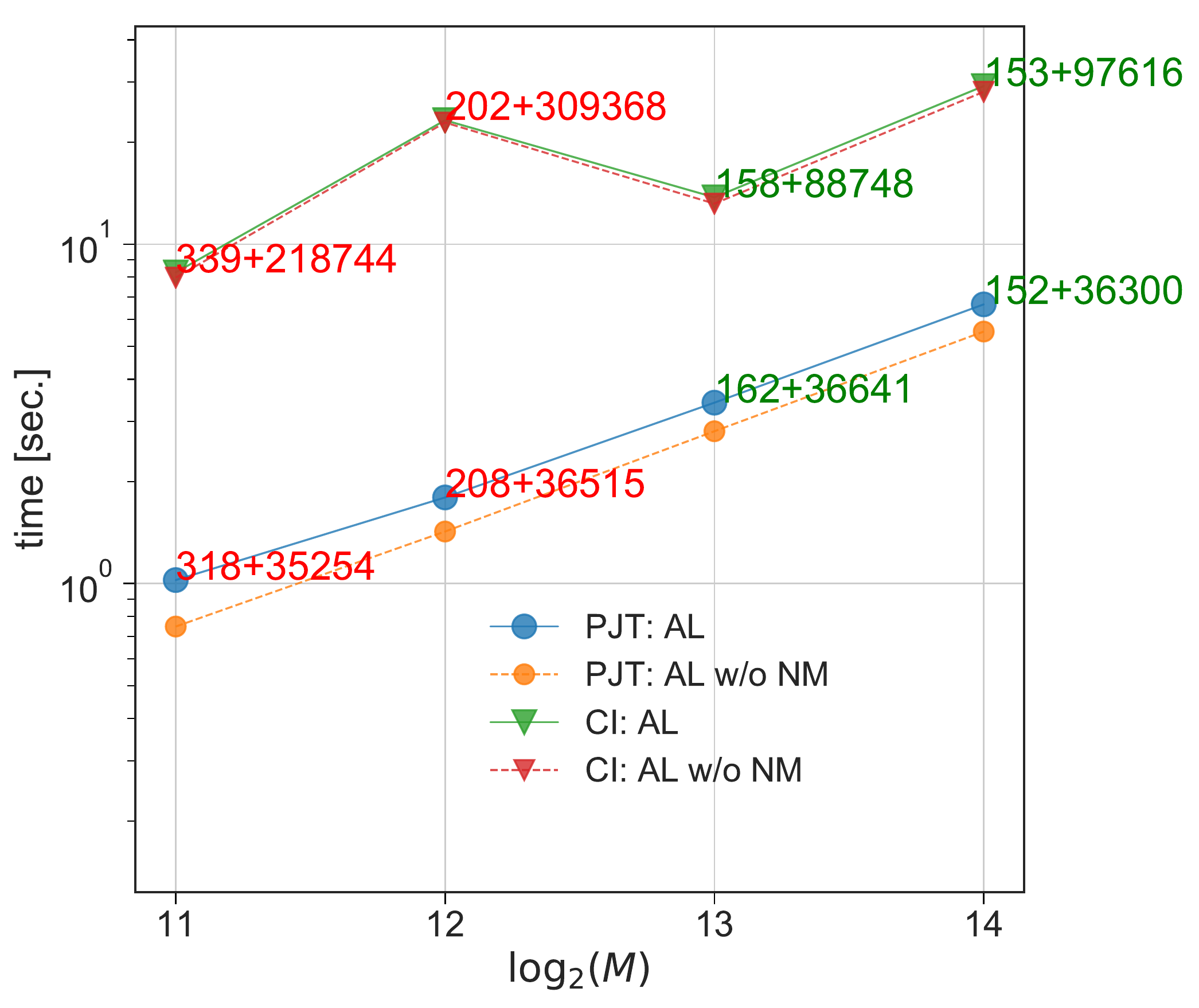}
 \includegraphics[width=0.445\linewidth]{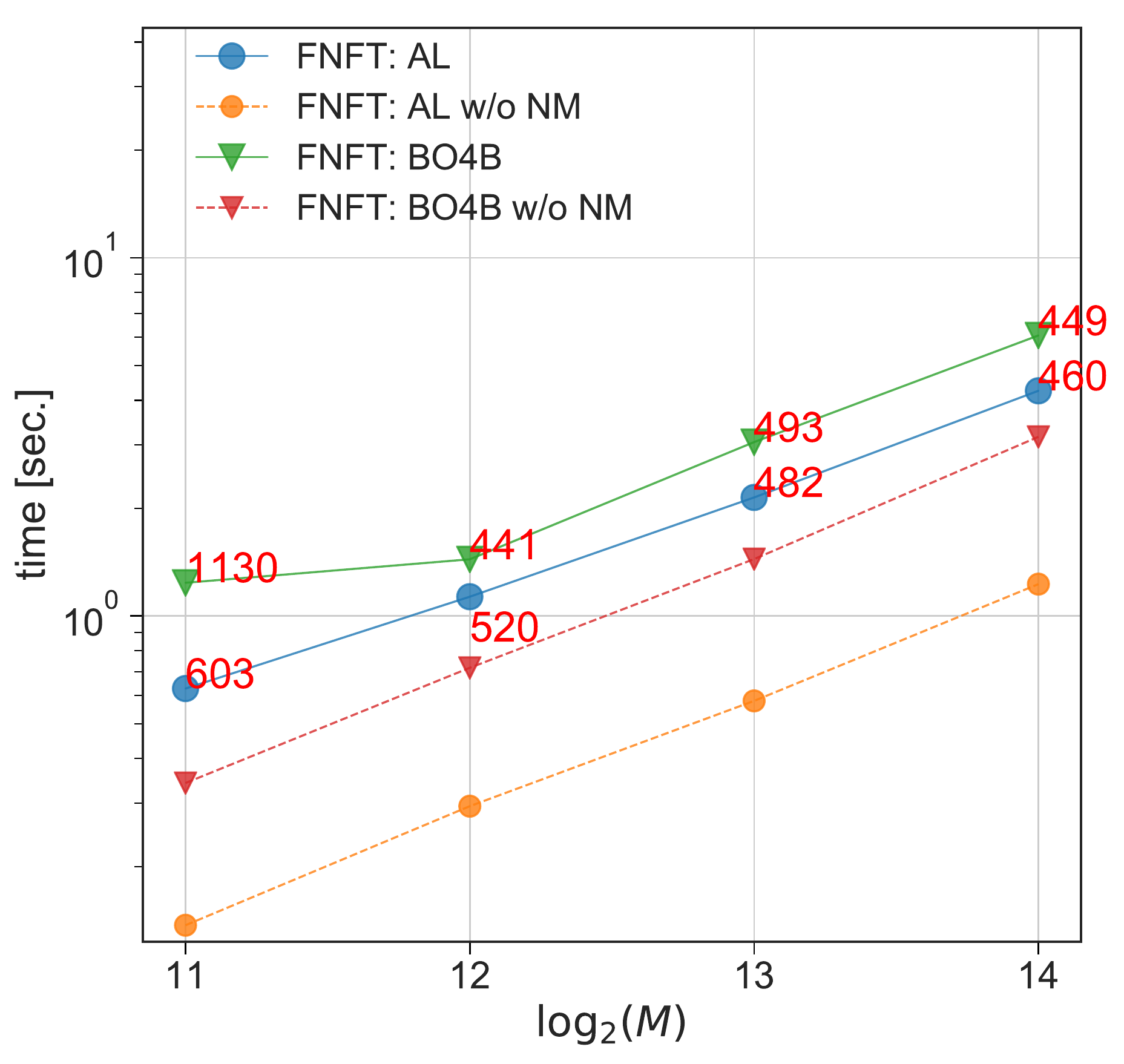}
 \vspace{-0.1cm}
 \caption{The runtime of various algorithms on the example with the 32-soliton solution. Solid lines indicate the total operating time, taking into account the iterative refinement of the discrete spectrum by the Newton method, dashed lines indicate the operating time without taking into account the iterative refinement. The left numbers represent the total number of calls to Newton's method. The right numbers on the plots indicate the total number of calls to the AL method for PJT and the CI method. The green colour of the numbers indicates that in this test case the discrete spectrum was correctly determined, and the red colour indicates the opposite.}
 \label{FIG:All_32}
\end{figure}

Figs.~\ref{FIG:All_20} and~\ref{FIG:All_32} show the results of comparing all the methods using the example of the oversoliton with parameters $ A = 20 $ and $ C = 0 $, as well as the 32-soliton solution. As the discrete spectrum size has increased, and the number of points in the signal has remained the same, the FNFT methods demonstrate a similar operating time with the PJT method. However, FNFT methods cannot correctly determine the discrete spectrum in any of the considered signal discretizations. The only exception is the oversoliton $ A = 20 $, $ C = 0 $, for which, for $ M = 2^{12} $, the FNFT method for the Boffetta-Osborne scheme managed to correctly find all the discrete eigenvalues, however, this is more likely an accident, since the initial approximations calculated by the FNFT method were far from true.
Generally speaking, in these examples, the AL method lacks accuracy for $ M < 2^{13} $ in order to correctly determine the discrete spectrum. Therefore, the correct answer was obtained by the PJT and CI methods only for $ M = 2^{13} $ and $ M = 2^{14} $. In other cases, for the CI method no more than 50,000 points were allocated to the initial contour $ \partial G $, which was still not enough.
In these examples, the PJT method turns out to be the undisputed leader, since it allows for $ M = 2^{14} $ to correctly determine the discrete spectrum of oversoliton with 20 discrete eigenvalues in 1 second, and the 32-soliton solution discrete spectrum in 6.6 seconds.

\section{Conclusion}

Summing up, we present a new method for finding the discrete spectrum of the direct Zakharov-Shabat problem, based on finding jumps in the argument of the coefficient $ a (\zeta) $. The method shows significant advantage over other methods when calculating a large discrete spectrum, both in speed and accuracy.
In many ways, the speed of the proposed algorithm is due to the optimized implementation of the AL method.

The method, of course, is not without drawbacks that can prove themselves in practical examples. But we are sure that the method has great prospects, therefore, by further refinement and development of the method, it will become a convenient tool for finding the discrete spectrum of the direct ZS problem. In addition, we plan to explore the possibility of generalizing it to various classes of analytic functions.

The method allows parallel implementation: firstly, the tracking of each of the trajectories can be performed independently of each other, and, secondly, each individual integration of the ZSP by the AL method can be parallelized.

Among the possible areas for further improvement of the method, one can note the use of FNFT to calculate the values of $ a (\zeta) $ in the upper complex half-plane, downsampling of the signal, and the use of the adaptive step to track the paths of the coefficient $ a (\zeta) $ jumps.

\section*{Funding}
Russian Science Foundation (RSF) (17-72-30006).

\bibliographystyle{unsrt}
\bibliography{library} 

\end{document}